\documentclass[12pt,a4paper]{article}
\usepackage{amsmath,amsfonts,latexsym}

\usepackage[all]{xy}
\def\ll{\:\raisebox{-1.6ex}{$\stackrel{\rm lt}{\rm t \to 0}$}\:}
\def\g{{\mathfrak{g}}}
\def\E{{\cal E}}

\def\PP{{\mathbb P}}

\title{Moduli of principal bundles in positive characteristic}
\author{ Vikram B. Mehta and S. Subramanian}
\date{}
\begin{document}
\maketitle

\noindent 1. Let $G$ be a semisimple algebraic group defined over an 
algebraically closed field $k$ of characteristic $p$.  In this article, we 
construct the moduli space of semistable principal $G$ bundles on a smooth
projective curve $X$ over $k$, of genus $g\geq 2$.  When the characteristic is 
zero, for example the field of complex numbers, these moduli spaces were first
constructed by A. Ramanathan (see [8], [9]).  Later on, different methods were 
used by Faltings (see [4]) and by Balaji and Seshadri (see [3]).  In these 
later papers, the emphasis was on proving the properness of the functor of 
semistable principal $G$-bundles. The method of Balaji and Seshadri was 
extended to positive characteristic by Balaji and Parameswaran (see [2]).\\

In this article, we follow A. Ramanathan's original approach.  The new input 
in the present  work is the use of low height and super low height 
representations (see [6]). These properties are used to ensure that the 
adjoint bundle of a semistable (respectively polystable)  $G$-bundle is a 
semistable  (respectively polystable) vector bundle.  However, the proof
that the moduli space is projective involves an  argument of lifting to
charateristic zero.  This is because one does not know if the kernel of the 
killing form of a Lie algebra over $k$ is a solvable ideal.  We try to
adhere to the notation of Ramanathan (see [8], [9]) as far as possible.

\medskip
\noindent
2. We recall the notion of a universal family and a universal space.

\paragraph*{Definition (2.1) (universal family)}: Let $ E \rightarrow X 
\times T$ be a family of semistable $G$-bunldes on $X$ parametrised by a 
scheme $T$.  Suppose an algebraic group $H$ acts on $T$ by $\alpha : H \times 
T \rightarrow T$ and also on $E$ as a group of $G$ bundle isomorphisms
compatible with $\alpha$.  We call $E \rightarrow X \times T$ a  universal 
family of semistable $G$-bundles with group $H$ if the following conditions 
hold:\\
(i) Given any family of semistable $G$ bundles $F \rightarrow X \times S$ and 
a point $s \in S$ there is an open neighborhood $U$ of $s \ \mbox{in}
 \ S$ and a morphism $ t: U \rightarrow T$ such that $ F \mid  X \times U$
is 
isomorphic to $t^* E  \rightarrow X \times U$ where $t^* E$ is the
pullback of $E$ by the morphism $ t : U \rightarrow T$.\\
(ii) Given two morphisms  $t_1, t_2:S \rightarrow T$ and an isomorphism
$\varphi : t^*_1 E \rightarrow t^*_2 E$  of $G$ bundles on $X \times  S$,
there is a unique morphism $h : S \rightarrow  H$  such that $t_2 = h t_1$
and $\varphi = (h \times t_1)^*(\alpha)$.

\paragraph*{Definition (2.2):}  Let $Sch$ denote the category of schemes 
over $k$, and let {\it Sets }denote the category of sets.  Let  
$\stackrel{\sim}{F}$ be the sheaf associated to the functor $F: Sch  
\rightarrow$ {\it Sets} which associates to a scheme  $T$ the set  of
isomorphism classes of semistable $G$ bundles parametrised by $T$. 
On morphisms $F$ is defined by pulling back.  Let $M$ be a scheme and $H$
an algebraic group acting on $M$ by $\alpha: H \times M \rightarrow M$.  Let 
$M/H$ be the sheaf associated  to the presheaf which assigns to a
$k$-scheme  $T$ the quotient set Hom $(T, M)/ \ \mbox{Hom} (T,H)$.  We say  
$M$ is a universal space with group $H$ if there is an isomorphism of sheaves
$\stackrel{\sim}{F} \rightarrow M/H$. With these definitions, we have the
following proportion.

\paragraph*{Proposition (2.3)}:  Suppose there is a universal space $M$
with  group $H$ for families of semistable $G$ bundles. If a good quotient of 
$M$ by $H$ exists, then it gives a coarse moduli scheme of semistable $G$-
bundles.

\paragraph*{Proof}:  The proof is the same as that of Proposition (4.5) in [9],
 using Lemma (4.3) and (4.4) of our article instead of the Proposition
(3.24) of [8].   \\
\hfill{QED}

In view of the above proposition, it is enough to construct a universal space 
for families of semistable $G$ bundles on $X$, and show that this universal 
space has a good quotient.  For this purpose, let $Z$ be the centre of $G$, 
and let $\g$ be the lie algebra of $G$.  We consider the sequence of group
homomorphisms

$$
G \rightarrow G/Z \rightarrow \ \mbox{Aut} \ (\g) \rightarrow GL(\g)
$$
where Aut$(\g)$ denotes the group of Lie algebra automorphisms of $\g$ and  
$GL (\g)$ is the group of linear automorphisms of the vector space $\g$.   
When the Killing form of $\g$ is non degenerate, $G/Z$ is the connected  
component of Aut $(\g)$(this follows from the fact that if the killing form 
of $\g$ is nondegerate, then every derivation of $\g$ is  inner, see [10], 
page 18).  A universal space $R_3$ for semistable $G$ bunldes is constructed 
in four steps.  We first need

\paragraph*{Definition (2.4)}:  Let $\alpha = \sum n_i \alpha_i$ be the 
highest root of the Lie algebra  $\g$, where $\alpha_i$ are simple roots and 
$n_i$ are positive integers. Then the height $h(G)$ of $\g$ is defined as 

$$
h(G) = 2 \sum n_i
$$

We have the following result from [6]:

\paragraph*{Theorem (2.5):}  Let $E \rightarrow X$ be a semistable $G$ bundle, 
and let $E(\g)$ be the adjoinst bundle of $E$.  If the characteristic $p$ is 
bigger than $h(G)$, then $E(\g)$ is a semistable vector bundle.  (see Theorem 
(2.2) in [6]).\\

Thus if $p > h(G)$, then every semistable $G$ bundle induces a semistable 
$GL(\g)$ bundle by the adjoint representation of $G$. We start with a  
universal family $R$ with group $GL(n)$ for semistable $GL(\g)$ bundles (which 
are simply semistable vector  bundles of rank = dim $\g$).  Since the
homomorphism Aut $(\g) \rightarrow GL(\g)$ is injective, we obtain a universal
family $R_1$ with group $GL(n)$ for semistable Aut $(\g)$ bundles.  This 
follows from the following lemma:

\paragraph*{Lemma (2.6)}: Let $A$ be a reductive group and $B$ a semisimple
group such that $B$ is a closed subgroup of $A$.  Then a universal family 
with group $H$ for families of semistable $A$ bundles gives a universal family
with group $H$ for families of semistable $B$ bundles, provided every
semistable $B$ bundle on $X$ induces a semistable $A$ bundle by the extension 
of structure group $B \rightarrow A$.

\paragraph*{Proof.} See Lemma (4.8.1) and Lemma (4.10) in [9]. \hfill{QED}

\paragraph*{Lemma (2.7):}The scheme $R_1$ is nonsingular and dim $R_1 = n^2 +
(r+1) (g-1)-g$ where $r = \dim G$, and $g$ is the genus of $X$.

\paragraph*{Proof:} See Lemma (4.13.4) in [9]. \hfill{QED}

The same procedure of applying Lemma (2.6) above to the injection $G/Z \subset 
\ \mbox{Aut} \ (\g)$, where $Z$ is the center of $G$ yields a universal family
$R_2$ with group $GL(n)$ for families of semistable $G/Z$ bundles on $X$.  We 
have

\paragraph*{Lemma (2.8):}  The natural morphism $\pi_2 : R_2 \rightarrow R_1$ 
is etale and finite, and hence $R_2$ is nonsingular.

\paragraph*{Proof:}  See Proposition (4.14) in [9]. \hfill{QED}

We now construct a universal space $R_3$ with group $GL(n)$ from $R_2$,
following an idea of Sols-Gomez (see [5]).  Let
$E_2 \rightarrow X \times R_2$ be the universal $G/Z$ bundle on
$X \times R_2$ and let $\rho: G \rightarrow G/Z$ be the natural projection.  
We asssume hereafter that the center $Z$ of $G$ is reduced, which is
the case if the characteristic $p$ of the field $k$ is bigger than $h(G)$.
Let $\Gamma (\rho, E_2)$ be the functor defined on the categoryy $Sch/R_2$
of schemes over $R_2$, which assigns to a scheme
$\varphi: S \rightarrow R_2$ the set of isomorphism classes of pairs
$(E, \psi)$ where $E$ is a principal $G$ bundle on $X \times S$ and $\psi$
is an isomorphism of the $G/Z$ bundle $E/Z$ with the $G/Z$ bundle
$\varphi^* E_2$.  Let $\stackrel{\sim}{\Gamma}(\rho, E_2)$ be the
sheafification of the functor $\Gamma(\rho, E_2)$.  We have a natural morphism 
$\pi_3: \stackrel{\sim}{\Gamma}(\rho, E_2) \rightarrow R_2$. We need a
definition.

\paragraph*{Definition (2.9):}  Let $\Gamma_1$ and $\Gamma_2$ be two
contravariant functors from {\it Sch} to {\it Sets} , and let $H$ be a 
finite reduced group, and let $\pi: \Gamma_1 \rightarrow \Gamma_2$ be 
a morphism of functors. Let $\stackrel{\sim}{\pi} : \stackrel{\sim}{\Gamma_1}
\rightarrow \stackrel{\sim}{\Gamma_2}$ be the morphism sheaves associated
 to $\pi$.  Then we say that $\stackrel{\sim}{\Gamma_1}$ is a principal
 $H$-bundle over $\stackrel{\sim}{\Gamma_2}$ if \\
(i) $H$ acts on $\Gamma_1$ and hence on $\stackrel{\sim}{\Gamma_1}$\\
(ii) the natural map

$$
\stackrel{\sim}{\Gamma_1} \times H \rightarrow \stackrel{\sim}{\Gamma_1}  
\displaystyle{\times_{\stackrel{\sim}{\Gamma_2}}} \stackrel{\sim}{\Gamma_1}$$ 
is an isomorphism of sheaves. We now have

\paragraph*{Proposition (2.10):} The natural map $\pi_3 : \stackrel{\sim} 
{\Gamma} (\rho,E_2)\rightarrow R_2$ is a principal $H^1 (X, Z)$ bundle on 
$R_2$, where  $H^1(X, Z)$ denotes the etale cohomology of $X$ with
coefficients in $Z$.

\paragraph*{Proof:}  One sees easily that there is a natural action of 
$H^1(X,Z)$ on $\Gamma(\rho,E_2)$.  We have to show that the action on 
$\stackrel{\sim}{\Gamma}(\rho,E_2)$ is free and the quotient is $R_2$.  We 
do this by identifying the fibre of $\pi_3: 
\stackrel{\sim}{\Gamma}(\rho,E_2) \rightarrow R_2$ with $H^1 (X, Z)$.
Let $S$ be a closed point of $R_2$ and $p : F \rightarrow X$ be corresponding
$G/Z$ bundle on $X$. $\stackrel{\sim}{\Gamma}(\rho, E_2)(S)$ now is the set 
of all pairs $(E, \alpha)$ such that $E \rightarrow X$ is a principal $G$ 
bundle on $X$ and $\alpha$ is an isomorphism $\alpha : E/Z \rightarrow F$. 
 We now check that $H^1(X, Z)$ acts freely and transitively on 
$\stackrel{\sim}{\Gamma} (\rho, E_2)(S)$.  The transitivity of the action is
easy to see.  To show that the action is free, we take a pair $(E, \alpha)$ in
$\stackrel{\sim}{\Gamma}(\rho, E_2)(S)$ and an element $t \in H^1(X, Z)$, and
show that $t(E, \alpha) \neq (E, \alpha)$.  If the $G$-bundle $t E$ is not 
isomorphic to $E$ we are through.  If $ t E$ is isomorphic to
$E$ as a $G$-bundle, let $\varphi: E \rightarrow t E$ be an isomorphism.  
Let $B$ denote a Borel subgroup of $G$ and $T$ denote a maximal torus of $G$
contained in $B$.  Since Ant$^0(G/B) \cong G/Z$ (this follows from the  
nondegeneracy of the Killing form of $\g$), we regard the $G/B$ fibration
$E/B \rightarrow X$ as equivalent to a $G/Z$ bundle.  Since $ E
\rightarrow  E/B$ is a $B$-bundle the canonical homomorphism $B \rightarrow 
T$ induces a $T$ bundle $E_T \rightarrow E/B$ on $E/B$.Now  $\alpha$ gives
an  
isomorphism $E/B \rightarrow F/B'$ of $G/B$ fibrations, where $B'$ is $B/Z, \ 
t \alpha$ gives an isomorphism $t E/B \rightarrow F/B'$ and $\varphi$ induces 
an isomorphism $\overline{\varphi}: E/B \rightarrow t E/B$.  We have to show 
that the diagram

$$
\xymatrix{
E/B \ar[d]_{\overline{\varphi} \cong} \ar[r]^\alpha_\cong & F/B\\
t E/B \ar[ru]^\cong_{t\alpha} &}
$$
is not commutative.  We note that there are two isomorphisms $ E/B \rightarrow
t E/B$, namely $\overline{\varphi}$ and $(t \alpha)^{-1} o \alpha$.  To show 
that the above diagram is not commutative, it is enough to check that the $T$-
bundle $E_T$ on $E/B$ pulls back, under $\overline{\varphi}$ and $(t \alpha)^{-1} o \alpha$, to non 
isomorphic $T$ bundles on $ t E/B$.  But this is clear as one of them is a 
translate of the other by $t$ and the action of $H^1(X,Z)$ on $H^1(E/B,
T)$ is 
fixed point free. \hfill{QED}

\paragraph*{Lemma (2.11):}  $\stackrel{\sim}{\Gamma}(\rho,E_2)$ is 
representable by a scheme $R_3$.

\paragraph*{Proof:}  Since the mapping of sheaves $\pi_3: \stackrel{\sim}
{\Gamma}(\rho, E_2) \rightarrow R_2$ is a principal $H^1(X, Z)$ bundle, it is 
easily seen that $\stackrel{\sim}{\Gamma}(\rho, E_2)$ is an algebraic space in 
the sense of Artin (see [1]).  If further follows from loc. cit. (see Theorem (3.3) in [1]), that this algebraic space is in 
fact a scheme. \hfill{QED}

\paragraph*{Remark (2.12):}  We observe that by contruction $\pi_3: R_3 
\rightarrow R_2$ is finite etale with structure group $H^1(X, Z)$.

\medskip 
\noindent
3. In order to construct the moduli space of semistable $G$ bundles on $X$,  
we have to construct a good quotient of $R_3$ by $GL(n)$.  In view of Lemma 
(2.8) and Remark (2.12) above,  it is enough to show that $R_1$ has a good 
quotient by appealing to Proposition (3.12) in [7].
Let $G_{n,r}$ be the Grassmannian of $r$ dimensional quotients of $n$- 
dimensional affine space over $k$,  where $r =\dim G = \dim \g$.  Let $Y= GL
(\g)/G_m \times \ \mbox{Aut} \ (\g)$ and let $Q$ be the universal quotient 
bundle on $G_{n,r}$.  Let $Q(Y)$ be the fibre bundle on $G_{n,r}$ with fibre
$Y$ associated to $Q$.  Let $\overline{Y}$ be the closure of $Y$ in the 
projective space $\PP (\g^* \otimes \g^* \otimes \g)$.  By the usual diagonal 
argument, we can choose an $N$-tuple $(x, \ldots, x_N)$ of points in $X$ for 
$N$ sufficiently large such that the evaluation morphism

$$
R_1 \rightarrow Q(\overline{Y})^N = \underbrace{Q(\overline{Y}) \times \cdots
 \times Q (\overline{Y})}_{N \ \mbox{factors}}
$$
is injective (see section (5.5) in [9]).  We can choose a polarisation (see 
(5.5) in [9]) so that if $Q(\overline{Y})^N_{ss}$ denotes the set of 
semistable points in the sense of geometric invariant theory for this 
polarisation then  $R_1$ maps into $Q(\overline{Y})^N_{ss}$ (see Lemma
(5.5.3) in [9]).

\paragraph*{Lemma (3.1):}  The map $R_1 \rightarrow Q(\overline{Y})^N_{ss}$ is 
proper.

\paragraph*{Proof:}  We check the valuative criterion for properness.  Let 
$k[[t]]$ denote the ring of formal power series in $t$ over $k$, and let 
$k((t))$ be its quotient field. We begin with a semistable vector bundle $V$ of
degree zero on $X \displaystyle{\times_{k}} \ \mbox{Spec} \ k ((t))$, which is 
also a bundle of Lie algebras, with Lie algebra structure isomorphic to that of
$\g$.  Let $ C \in \PP (\g^* \otimes \g^* \otimes \g)$ be the structure
tensor of the Lie algebra $\g $.  Then the Lie algebra structure on $V$ has 
structure  tensor  $M((t)) C$ where $M((t))$ is an element of $GL(r, k
((t))$.  If $\overline{V}$ is any semistable extension of $V$ to $X 
\displaystyle{\times_{k}} \ \mbox{Spec}k[[t]]$, then the special fibre $V_o$ 
of $\overline{V}$ is also a bundle of Lie algebras with structure tensor
$\ll M((t)) C$ (the limit is taken in $ \PP(\g^* \otimes \g^* \otimes \g))$.  
We now show that this limit Lie algebra structure is isomorphic  to that 
defined by $C$ (namely,  $\g$ itself).
Let $W$ be the ring of Witt vectors of $k$, let $K$ be the quotient field of 
$W$.  The Lie algebra $\g$ has a form $\g_W$ over $W$,  whose structure 
tensor we denote by $C_W \in \PP(\g^*_W \otimes \g_W^* \otimes \g_W)$.  We now 
lift the matrix $M((t))$ to a matrix $M_W((t))$ which is an element of
$GL(r, W((t))$.  We therefore obtain a Lie algebra $\g_{W,t}$ with structure 
tensor $M_W((t)) C_W$ over $W((t))$. The group scheme Aut $(\g_{W,t})$ of Lie algebra automorphisms of $\g_{W,t}$ over $W((t))$ is 
smooth over $W((t))$ (see the proof of the Theorem in (11.5) in [2]).  Let 
$\g_t$ be the Lie algebra $\g_{W,t}$ mod $p$ over $k((t))$.  Then by Lemma 1, 
p. 175 of [11], the principal Aut( $\g_t)$ bundle on $X \times \
\mbox{spec}\ 
k((t))$ can be lifted to a principal Aut $(\g_{W,t})$ bundle $E_W$ on $X_W 
\displaystyle{\times_{W}} W((t))$ where $X_W$ is a form of the curve $X$ over
the Witt ring $W$.  Since $V$ is semistable on $X \times \ \mbox{Spec}\ k((t))$
the generic fibre $E_K$ of $E_W$ is semistable on $X_K\times K((t))$.  Now by 
the key lemma (5.6) in [9], $V_K$ has a semistable extension
$\overline{V}_K$ to $X_K \times K [[t]]$, where $V_K$ is the vector bundle 
associated to $E_K$, such that the Lie algebra structure on the special fibre 
$\overline{V}_{K,0}$ at $t=0$ is isomorphic to the Lie algebra structure on 
$V_K$.  Hence we obtain that $\ll M_W((t)) C_W$ gives a Lie algebra isomorphic 
to that given by $C_W$.  It follows by reducing mod $p$ that  the Lie
algebra structure $ \ll M((t)) C$ is isomorphic to the Lie algebra structure 
$C$. Hence $V_0$ is also a principal Aut $(\g)$ bundle (semistable by 
construction).  This shows that the map $R_1 \rightarrow Q (\overline{Y})
^N_{ss}$ is proper in positive characteristic. \hfill{QED}

After the above lemma, we see that a good quotient of $R_1$ exists and hence 
as remarked before, a good quotient of $R_3$ exists.  This completes the 
construction of the moduli space of semistable $G$ bundles.  Further, since 
the quotient of $R_1$ is projective (by the Lemma above ) and the map 
$R_3 \rightarrow R_1$ is finite and $GL(n)$ equivariant it follows that the 
quotient of $R_3$,  namely the moduli space,is also projective. It only 
remains to identify the closed points of the moduli space.

\medskip 
\noindent
4.  A description of the closed points of the moduli space is given by

\paragraph*{Theorem (4.1):}  The closed points of the moduli space of 
semistable $G$ bundles are isomorphism classes of polystable $G$ bundles.

\paragraph*{Proof:}  The proof follows from the following two lemmas. 
\hfill{QED}

Before the lemmas, we recall the following definition from [8].

\paragraph*{Definition (4.2):}  Let $E \rightarrow X$ be a semistable $G$ 
bundle on $X$.  If $E$ is stable, we define gr $E=E$.  If $E$ is not stable, 
there is an  admissible reduction of structure group $E_P\subset E$ to a 
parabolic $P$,  such that the  bundle $E_M$ induced from $E_P$ by the
projection $P \rightarrow M$ ($M$ being a Levi of $P$) is stable, and gr $E$ 
is the $G$-bundle obtained from $E_M$ by the inclusion $M \subset G$ (see 
(3.12) in [8]).

We now have

\paragraph*{Lemma(4.3):}  Let $ E \rightarrow X$ be a semistable $G$ bundle on 
$X$ regarded as a point of $R_3$.  Let $O(E)$ denote the orbit of $E$ under 
the group $GL(n)$ acting on $R_3$, and let $\overline{O(E)}$ denote the 
closure of the orbit in $R_3$. Then gr $E \in \overline{O(E)}$.

\paragraph*{Proof:}  The proof is by constructing a family of semistable $G$ 
bundles $\E \rightarrow X \times A^1$, where $A^1$ is the affine line over 
$k$, such that $\E_t$ is isomorphic to $E$ if $ t \neq 0$ in $A^1$, and 
$\E_0$ is isomorphic to gr $E$, where $\E_t = \E \mid X \times t, \ t \in
A^1$. This is done as in the proof of Proposition (3.24) (ii) in [8].  
\hfill{QED}

\paragraph*{Lemma (4.4):}  Let $ E\rightarrow X$ be a semistable $G$ bundle 
regarded as a point of $R_3$.  Then the orbit $O(\mbox{gr} \ E)$ of gr $E$ is
closed in $R_3$ if $ p > r h (G)$, where $ r = \dim G, \ p=$  characteristic 
of $k, \ h(G)$ is the height as in Definition (2.4) above

\paragraph*{Proof:}  It is enough to prove $(*)$ if gr $F \in
\overline{O \ 
(\mbox{gr} \ E)}$ then gr $F \cong \ \mbox{gr} \ E$. We take the
adjoint
representation of $G$ and consider the associated vector bundles $(\mbox{gr}F)
(\g)$ and $(\mbox{gr}E) (\g)$.  The hypothesis  on $p$ ensures that 
$(\mbox{gr} F) (\g)$ and $(\mbox{gr} E) (\g)$ are polystable vector bundles 
(see Theorem (4.2) in [6]).  If follows from complement (5.8.1) in [7],
that $\mbox{gr} F(\g)$ is isomorphic to $(\mbox{gr} E) (\g)$ as vector 
bundles. Since $G/Z \subset GL (\g)$ is a faithful representation by the 
adjoint action, it follows that gr $(E/Z)$ is isomorphic to  gr$(F/Z)$ as
$G/Z$ bundles.  Hence gr $(F/Z)$ belongs to the $GL(n)$ orbit of gr
$(E/Z)$ in
$R_2$.  However, since $\pi_3 : R_3 \rightarrow R_2 $ is a finite $GL(n)$ 
equivariant map, it follows  that gr $F$ belongs to the $GL(n)$ orbit of gr
 $E$, i.e.,  gr $F$ is isomorphic to gr $E$ as $G$ bundles.  This
completes the proof of $(*)$ and hence of the lemma. \hfill{QED}

\medskip
\noindent
5.  To summarise, we need the semisimple group $G$ to satisfy the following 
conditions:\\
(i) the centre $Z$ is reduced.\\
(ii)  the Killing form is non degenerate\\
(iii)  $p > 2h(G) + 3$.\\
(iv)  $p > \dim G. h (G)$

We remark that if $ p > \dim G. h (G)$ then all the other conditions are 
satisfied.  Hence a projective moduli space of semistable $G$ bundles exists 
if  $p > \dim G. h(G)$, and the closed points of this moduli space are
isomorphism closes of polystable bundles.

\newpage
\section*{References}

\begin{enumerate}
\item M. Artin, The Implicit Function theorem in algebraic geometry, Algebraic
geomety, Proceedings of the Bombay Colloquium 1968, Oxford University Press.
\item V. Balaji, A.J. Parameswaran, Semistable Principal Bundles-II, to appear.
\item V. Balaji, C.S. Seshadri, Semistable Principal Bundles -I, to appear.
\item G. Faltings, Stable $G$-bundles and Projective connections, J. Algebraic
Geometry, 2 (1993), p.507-568.
\item T.L. Gomez, I. Sols, Moduli space of Principal sheaves over projective 
varieties, preprint.
\item V.B. Mehta, A.J. Parameswaran, Geometry of low height representations, 
Proceedings of the International Colloquium on Algebra, Arithmetic and 
Geometry, Mumbai 2000, p. 417-426.
\item P.E. Newstead, Lectures on Introduction to Moduli Problems and orbit 
spaces, TIFR Bombay 1978, Narosa Publishing House.
\item  A. Ramanathan, Moduli of Principal Bundles over Algebraic Curves -I, 
Proc. Indian Acad. Sci (Math.Sci). Vol 106, No.3, August 1996, pp. 301-328.
\item A. Ramanathan, Moduli of principal bundles over algebraic curves-II, 
Proc. Indian Acad. Sci (Math. Sci), Vol 106, No.4, November 1996, pp. 421-449.
\item G.B. Seligman, Modular Lie algebras, Springer-Verlag Berlin Heidelberg 
New York 1967.
\item C.S. Seshadri, Desingularisation of the moduli varieties of vector 
bundles on curves, Proceedings of the International Symposium on algebraic 
geometry, Kyoto, January 1977, Kinokuniya Book Co. 1978, pp 155-184.
\item W.H. Hesselink, Uniform instability in reductive groups, Journal fur die 
reine und ang. Math. 303-304, 1978, pp. 74-96.\\[10mm]
\end{enumerate}

\noindent
vikram@math.tifr.res.in\\[-1mm]
subramnn@math.tifr.res.in\\[-1mm]
School of Mathematics\\[-1mm]
Tata Institute of Fundamental Research\\[-1mm]
Homi Bhabha Road\\[-1mm]
Mumbai 400 005 India

\end{document}